\documentclass[a4paper,twoside]{article}

\usepackage{lmodern}
\usepackage[T1]{fontenc}
\usepackage[ansinew]{inputenc}
\usepackage{latexsym}
\usepackage{amssymb}
\usepackage{amscd}
\usepackage{dsfont}
\usepackage{amsthm}
\usepackage{enumerate}
\usepackage[intlimits]{amsmath} %
\usepackage[all,arc,curve,matrix]{xy}
\usepackage{hyperref}
\usepackage{fancyhdr}

\makeatletter
\def\section{\@startsection{section}{1}{\z@}%
    {-21dd plus-4pt minus-4pt}{10.5dd plus 4pt
     minus4pt}{\normalsize\bfseries\boldmath}}
\def\subsection{\@startsection{subsection}{2}{\z@}%
    {-21dd plus-4pt minus-4pt}{10.5dd plus 4pt
     minus4pt}{\normalsize\itshape}}
\def\subsubsection{\@startsection{subsubsection}{3}{\z@}%
    {-13dd plus-4pt minus-4pt}{-5.5pt}{\normalsize\itshape}}
\def\paragraph{\@startsection{paragraph}{4}{\z@}%
    {-13pt plus-4pt minus-4pt}{-5.5pt}{\normalsize\itshape}}
\makeatother

\newcommand{\sect}{\section}

\newtheorem{thm}{Theorem}[section]
\newtheorem{prop}[thm]{Proposition}
\newtheorem{lem}[thm]{Lemma}
\newtheorem{cor}[thm]{Corollary}

\theoremstyle{definition}

\newtheorem{defn}[thm]{Definition}
\newtheorem{ex}[thm]{Example}

\newcommand{\email}[1]{\href{mailto:#1}{#1}}
\newenvironment{acknowledgement}{\small \setlength{\parindent}{0pt}\emph{Acknowledgements.}}{}

\newenvironment{myproof}{\begin{proof}}{\end{proof}}

\newcommand{\commdiag}[2]
{\hspace{8mm}\begin{xy}
  \xymatrix#1{#2}\end{xy}}

\renewcommand{\bar}{\overline}
\newcommand{\Qbar}{{\bar{\Q}}}

\newcommand{\N}{\ensuremath{\mathds{N}}}
\newcommand{\R}{\ensuremath{\mathds{R}}}
\newcommand{\Q}{\ensuremath{\mathds{Q}}}
\newcommand{\Z}{\ensuremath{\mathds{Z}}}
\newcommand{\C}{\ensuremath{\mathds{C}}}
\newcommand{\F}{\ensuremath{\mathds{F}}}
\renewcommand{\H}{\ensuremath{\mathds{H}}}
\renewcommand{\P}{\ensuremath{\mathds{P}}}
\renewcommand{\L}{\ensuremath{\mathcal{L}}}

\newcommand{\B}{\ensuremath{\mathcal{B}}}

\newcommand{\A}{\ensuremath{\mathcal{A}}}
\newcommand{\G}{\ensuremath{\mathcal{G}}}
\newcommand{\T}{\ensuremath{\mathcal{T}}}
\newcommand{\nt}{\ensuremath{^{\mathrm{nt}}}}
\newcommand{\nc}{\ensuremath{^{\mathrm{nc}}}}

\newcommand{\Mcal}{\ensuremath{\mathcal{M}}}
\newcommand{\Hcal}{\ensuremath{\mathcal{H}}}
\newcommand{\Ncal}{\ensuremath{\mathcal{N}}}
\newcommand{\Ccal}{\ensuremath{\mathcal{C}}}

\DeclareMathOperator{\PGL}{PGL}
\DeclareMathOperator{\PSL}{PSL}

\DeclareMathOperator{\SL}{SL}

\DeclareMathOperator{\ord}{ord}
\DeclareMathOperator{\Aut}{Aut}

\DeclareMathOperator{\Gal}{Gal}

\DeclareMathOperator{\id}{id}

\DeclareMathOperator{\Epi}{Epi}

\DeclareMathOperator{\Inn}{Inn}
\DeclareMathOperator{\Out}{Out}

\DeclareMathOperator{\Stab}{Stab}

\newcommand{\RA}{\Rightarrow}

\newcommand{\onto}{\twoheadrightarrow}
\newcommand{\x}{\cdot}
\newcommand{\restrict}[1]{|_{#1}}

\newcommand{\NT}{\vartriangleleft}

\newcommand{\set}[1]{\left\{#1\right\}}
\newcommand{\abs}[1]{\left|#1\right|}

\newcommand{\defi}[1]{\emph{#1}\index{#1}}
\newcommand{\defindex}[2]{\emph{#1}\index{#2}}
\newcommand{\zitat}[2]{\cite{#2}, #1}
\newcommand{\mat}[1]{\begin{pmatrix}#1\end{pmatrix}}
\newcommand{\smat}[1]{\left(\begin{smallmatrix}#1\end{smallmatrix}\right)}

\newcommand{\gen}[1]{\left<#1\right>}
\newcommand{\mattwo}[4]{\left(\begin{matrix}#1 & #2 \\ #3 & #4\end{matrix}\right)}

\newcommand{\OriSquareBorderMark}[6]
{
\begin{xy}
(0,0)="Pos";
"Pos"+(1,0) **@{#6}; ?*h!C!/_0pt/{\text{\tiny #5}}, 
"Pos"+(1,1) **@{#6}; ?*h!C!/_0pt/{\text{\tiny #2}},
"Pos"+(0,1) **@{#6}; ?*h!C!/_0pt/{\text{\tiny #3}},
"Pos" **@{#6}; ?*h!C!/_0pt/{\text{\tiny #4}},
"Pos"+(0.5,0.5) *h{\text{#1}};
\end{xy}
}

\newcommand{\KappesOrigami}
{
\begin{minipage}{4cm}
\vspace{1mm}
\begin{xy}
<0.8cm,0cm>:
(0,0)*{\OriSquareBorderMark{}{}{|}{=}{|}{-}};
(1,0)*{\OriSquareBorderMark{}{}{}{}{||}{-}};
(1,1)*{\OriSquareBorderMark{}{}{||}{$-$}{}{-}};
(2,0)*{\OriSquareBorderMark{}{=}{}{}{|||}{-}};
(2,1)*{\OriSquareBorderMark{}{}{|||}{}{}{-}};
(3,1)*{\OriSquareBorderMark{}{$-$}{||||}{}{||||}{-}};
\end{xy}
\vspace{1mm}
\end{minipage}
}

\newcommand{\GraphOrigami}[3]
{\begin{minipage}{1cm}\vspace{-5mm}\begin{graph}
[]\bullet[]
(-_{#2}@(ul,dl)[],
-@0_{#1\;}@(r,d)[]		
:^{#3}[r]
)
\end{graph}\end{minipage}}

\newcommand{\GraphOrigamiAut}[7]
{\begin{minipage}{1cm}\begin{graph}
[]\bullet[]
(
:^{#4}[ul],
:_{#5}[dl]
)
-@0_{#1}@(r,d)[]		
-^{#3}[r]		
-@0_{#2}@(d,l)[]		
[]\bullet[]
(
:_{#6}[ur],
:^{#7}[dr]
)
\end{graph}\end{minipage}}

\newcommand{\GraphTriangle}
{\begin{minipage}{1cm}\begin{graph}
[]\bullet[]
:_{C_{3}}[l]
[r]\bullet[]
-@0_{A_5}@(r,d)[]		
-^{D_5}[r]
[]\bullet[]		
-@0_{\;\;D_{5}\!\!\!}@(r,d)[]		
:_{C_{5}}[d]
[u]
-^{D_5}[r]		
[]\bullet[]
-@0_{A_5}@(r,d)[]		
:^{C_{3}}[r]
\end{graph}\end{minipage}}

\fancypagestyle{title}{
\fancyhf{}
\fancyhead[L]{To appear in manuscripta mathematica \hfill DOI: \href{http://dx.doi.org/10.1007/s00229-010-0379-8}{10.1007/s00229-010-0379-8}\vspace{1pt}\\
The final publication is available at \href{http://www.springerlink.com/content/t5r648x630230437/fulltext.pdf}{www.springerlink.com}
}
\fancyfoot[L]{\footnotesize 
K. Kremer:    Institute for Algebra and Geometry,
               Karlsruhe Institute of Technology (KIT),
              76128 Karlsruhe, Germany.
              e-mail: \email{kremer@kit.edu} \medskip 
}

\setlength{\footskip}{24pt} 
}

\fancypagestyle{plain}{
\fancyhead[LE,RO]{\small \thepage}
\fancyhead[RE]{\small Karsten Kremer}
\fancyhead[LO]{\small Normal origamis of Mumford curves}
\fancyfoot{}

}

\begin{document}
\thispagestyle{title}

{\setlength{\parindent}{0pt}\phantom{x}\vspace{0.5cm}

Karsten Kremer \bigskip

{\Large \bf Normal origamis of Mumford curves} \vspace{10pt}

{\small 19 February 2010} \bigskip

{\small
{\bf Abstract.}
An origami (also known as square-tiled surface) is a Riemann surface covering a torus with at most one branch point. Lifting two generators of the fundamental group of the punctured torus decomposes the surface into finitely many unit squares. By varying the complex structure of the torus one obtains easily accessible examples of Teichm\"uller curves in the moduli space of Riemann surfaces.
The $p$-adic analogues of Riemann surfaces are Mumford curves. A $p$-adic origami is defined as a covering of Mumford curves with at most one branch point, where the bottom curve has genus one. A classification of all normal non-trivial $p$-adic origamis is presented and used to calculate some invariants. These can be used to describe $p$-adic origamis in terms of glueing squares.\\
}
}\medskip

\hrule\medskip

\sect{Introduction}
An origami is a covering of a nonsingular complex projective curve over an elliptic curve which may be ramified only over zero. 
We will see in Section \ref{sct:complex} that such an object defines a curve (called origami-curve) in the moduli space of Riemann surfaces. Origamis have been studied recently for example by Lochak \cite{Loch}, Zorich \cite{Zor}, Schmith\"usen \cite{Schm2} and Herrlich \cite{HS1}. Several equivalent definitions for origamis are in use, see \zitat{Ch.~1}{KreDiss}. Our definition can be generalized to other ground fields, such as the $p$-adic field $\C_p$. In the complex world we get every Riemann surface as a quotient of an open subset $\Omega$ of $\P^1(\C)$ by a discrete subgroup $G$ of $\PSL_2(\C)$. In the $p$-adic world the analogues of Riemann surfaces, which admit a similar uniformization $\Omega/G$, are called \emph{Mumford curves}. But contrary to the complex world not every nonsingular projective curve over $\C_p$ is a Mumford curve. Mumford curves have been thoroughly studied; two textbooks giving a comprehensive introduction are \cite{GvdP} and \cite{FvdP}.
\pagebreak[2]

As Mumford curves are the $p$-adic analogues of Riemann surfaces we define $p$-adic origamis to be coverings of Mumford curves with only one branch point, where the bottom curve has genus one. In Section \ref{sct:p-adic_origamis} we classify all normal non-trivial $p$-adic origamis. This is done using the description of the bottom curve as an orbifold $\Omega/G$, where $\Omega \subset \P^1(\C_p)$ is open and $G$ is a group acting discontinuously on $\Omega$. These groups and the corresponding orbifolds can be studied by looking at the action of $G$ on the Bruhat-Tits tree of a suitable subfield of $\C_p$ and the resulting quotient graph of groups. This has been done by Herrlich \cite{Her}, and more recently by Kato \cite{Kato2} and Bradley \cite{Brad}.

We will see in Section \ref{sct:p-adic_origamis} that all normal $p$-adic origamis with a given Galois group $H$ are of the type $\Omega/\Gamma \to \Omega/G$ with the following possible choices for the groups $\Gamma$ and $G$: The quotient graph of $G$ can be contracted to
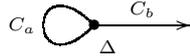
\begin{center}\vspace{4mm}
\GraphOrigami{\Delta}{C_a}{C_b}
\vspace{-2mm}
\end{center}
for a $p$-adic triangle group $\Delta$ (where the single vertex represents a subtree with fundamental group $\Delta$, the arrow indicates an end of the graph, and $C_a$ and $C_b$ are finite cyclic groups of order $a$ resp. $b$), which means that $G$ is isomorphic to the fundamental group of this graph, i.e.
\[G \cong \left<\Delta,\gamma;\, \gamma\alpha_1=\alpha_2\gamma\right> \text{ with $\alpha_i \in \Delta$ of order $a$.}\]
$\Gamma$ is the kernel of a morphism $\varphi : G \to H$ which is injective when restricted to the vertex groups of the quotient graph of $G$. The ramification index of the $p$-adic origami is then $b$. We have a similar result (Theorem \ref{thm:p-adic_origami_aut}) for the automorphism group of the $p$-adic origami. 

Given a $p$-adic origami which is defined over $\Qbar$ we can change the ground field to $\C$ and know that there our origami can be described as a surface glued from squares. Actually doing this is usually hard, because we would have to work out equations for the Mumford curves and for the complex curves corresponding to the Riemann surfaces. Nevertheless in many cases we can find out which complex origami-curve belongs to our $p$-adic origami, as mostly the curve is already uniquely determined by fixing the Galois group. In Section \ref{sct:unique_curve} we prove that this is true for the Galois groups $D_n \times \Z/m\Z$ and $A_4 \times \Z/m\Z$ (with $n,m \in \N$ and $n$ odd). 

We will also discuss some cases where this does not work, i.e. where there are several origami-curves of origamis with the same Galois group. For groups of order less than or equal to 250 this happens only for 30 groups.
To construct examples we can take an origami-curve which is not fixed by a certain element $\sigma$ of the absolute Galois group $\Gal(\Qbar/\Q)$. As the action of $\Gal(\Qbar/\Q)$ on origami-curves is faithful by \zitat{Theorem 5.4}{Mol} we can find such a curve for any given $\sigma$. In this case of course both the curve and its image contain origamis with the same Galois group, and we suspect that all other known invariants of origami-curves are equal as well.

\setlength{\textheight}{19.0cm}
\sect{Complex origamis}\label{sct:complex}

An \defi{origami} is a covering $p : X \to E$ of degree $d$ of a connected surface $X$ over the torus $E=\R^2/\Z^2$ which may be ramified only over $0 \in E$. We can lift the unit square defining $E$ to $X$. This yields a decomposition of $X$ into $d$ copies of the unit square.  Removing the ramification points leads to an unramified restriction $p : X^* := X \setminus p^{-1}(0) \to E^* := E \setminus \set{0}$ of degree $d$.
\pagebreak[4]

The covering is \defi{normal} if a (Galois) group $G$ acts on $X$ such that $E = X/G$.
The monodromy of such a covering is by definition the action\footnote{\label{foot:ori:right_to_left}Note that if we want to consider elements $\alpha,\beta \in \pi_1(E^*,P)$ as permutations of the fiber $p^{-1}(P)$, we need $\alpha\beta$ to be the path \emph{first} along $\beta$ and \emph{afterwards} along $\alpha$. This may not be an intuitive way to define multiplication in $\pi_1(E^*,P)$, but otherwise the group would not act on the fiber from the left.} of the fundamental group $\pi_1(E^*,\bar{P})$ on the fiber $p^{-1}(\bar{P})=G \x P$ over any basepoint $\bar{P}=p(P)$ with $P \in X^*$. Without loss of generality we can choose both coordinates of $\bar{P}$ to be non-zero in $\R/\Z$. The fundamental group $\pi_1(E^*,\bar{P})$ is isomorphic to the free group generated by $x$ and $y$, where $x$ is the closed path starting at $\bar{P}$ in horizontal direction, and $y$ is the closed path starting at $\bar{P}$ in vertical direction. Now we have $F_2$ acting on the orbit $G \x P$, and this action has to be compatible with the group action of $G$, thus the monodromy is defined by a homomorphism $f : F_2 \to G$. If we think of the squares making up $X$ labelled by elements of $G$ then the right neighbor of $g \in G$ is $f(x)g$, and its upper neighbor is $f(y)g$. As the surface is connected $f$ has to be surjective. Of course the monodromy map is only well-defined up to an automorphism of $G$.

We identify the torus $E$ with $\C / \Z[i]$, thus an origami $p : X \to E$ becomes a Riemann surface using the coordinate charts induced by $p$. In fact, we get a lot of Riemann surfaces: for every $A \in \SL_2(\R)$ we can define the lattice $\Lambda_A = A\x \Z^2$ and the homeomorphism $c_A : \R^2/\Z^2 \to \R^2/\Lambda_A =: E_A, x \mapsto A\x x$. The identification of $\R^2$ with $\C$ then leads to new coordinate charts induced by $p_A := c_A \circ p$. We get again a complex structure on the surface $X$ which we denote by $X_A$. 

If the torus $E_A$ is isomorphic to our torus $E = E_I$ as a Riemann surface, then $p_A : X \to E_A \cong E$ defines another origami. We thus get an action of $\SL_2(\Z)$ on the set of all origamis. If $\varphi \in \Aut(F_2)$ is a preimage of $A \in \SL_2(\Z) \cong \Out(F_2) = \Aut(F_2)/\Inn(F_2)$, then the action of $A$ on the monodromy $f : F_2 \to G$ of a normal origami is given by $f \circ \varphi^{-1}$ (\zitat{Prop. 1.6}{KreDiss}). 

The moduli space $\Mcal_{g,n}$ is the set of isomorphism classes of Riemann surfaces of genus $g$ with $n$ punctures, endowed with the structure of an algebraic variety. An origami defines a  subset $c(O) = \set{X_A : A \in \SL_2(\R)}$ in $\Mcal_{g,n}$. This set turns out to be an algebraic curve (\zitat{Prop. 3.2 ii)}{Loch}), called an origami-curve. The origami-curves of two origamis of the same degree $d$ are equal if and only if the two origamis are in the same $\SL_2(\Z)$ orbit (\zitat{Prop. 5 b)}{HS1}). The set of normal origamis with a given Galois group $G$ can thus be identified with the set $ \Aut(G) \backslash \Epi (F_2,G) / \Out(F_2)$.

We would like to define the automorphism group of an origami as an invariant of the origami-curve. Therefore we call a bijective map $\sigma : X \to X$ an \defi{automorphism} of the origami, if it induces for every $A \in \SL_2(\R)$ via $X \to E \to E_A$ a well-defined automorphism on $E_A$. An automorphism is called \defi{translation} if it induces the identity on $E$. An origami of degree $d$ is normal if and only if it has $d$ translations. In this case the group of translations is isomorphic to the Galois group $G$ (\zitat{Prop. 3.12}{KreDiss}).

\sect{Discontinuous groups}\label{sct:discontinuous}

After defining Mumford curves we will construct in this section the \emph{Bruhat-Tits-Tree} $\B$ for an extension of $\Q_p$ using a quite concrete definition from $\cite{Her}$. A Mumford curve is closely related to the quotient graph of an action of $G$ on a subtree of $\B$. Often one defines a suitable subtree such that the quotient becomes a finite graph, but instead we will follow Kato \cite{Kato2}, who uses a slightly larger quotient graph, which can be used to control the ramification behavior of coverings of Mumford curves.

\begin{defn} Let $k$ be a field which is complete with respect to a non-ar\-chi\-me\-de\-an valuation and $G$ a subgroup of $\PGL_2(k)$. 
A point $x \in \P^1(k)$ is called a \defi{limit point} of $G$, if there exist pairwise different $\gamma_n \in G$ $(n \in \N)$ and a point $y \in \P^1(k)$ satisfying $\lim \gamma_n(y) = x$. The set of limit points is denoted by $\L(G)$. 

$G$ is a \defi{discontinuous group}, if $\Omega(G) := \P^1(k) \setminus \L(G)$ is nonempty and for each $x \in \P^1(k)$ the closure $\overline{Gx}$ of its orbit is compact. A discontinuous group $G$ is called a \defi{Schottky group} if it is finitely generated and has no non-trivial elements of finite order. Every Schottky group is free (\zitat{Theorem I.3.1}{GvdP}).

A discontinuous group $G$ acts properly discontinuously on $\Omega(G)$. For a Schottky group $G$ we know from \zitat{Theorem III.2.2}{GvdP} that the quotient $\Omega(G)/G$ is the analytification of an algebraic curve. Such a curve is called a \defi{Mumford curve}. 
\end{defn}
If an arbitrary group $G \subset \PGL_2(k)$ contains a discontinuous group $G'$ of finite index, then $\L(G)=\L(G')$ and $G$ is also discontinuous. We know from \zitat{Ch. I, Theorem 3.1}{GvdP} that every finitely generated discontinuous group contains a Schottky group as a subgroup of finite index.

Let $k \subset \C_p$ be a finitely generated extension of $\Q_p$. Then the set of absolute values $\abs{k^\times} := \set{\abs{x}:x\in k^\times}$ is a discrete set in $\R^\times$. For $r \in \abs{k^\times}$ and $x \in k$ let $B(x,r) := \set{y \in k: \abs{x-y}\leq r}$ be the ``closed'' ball\footnote{Note that as $\abs{k^\times}$ is discrete the ball $B(x,r)$ is both open and closed for the topology induced by the $p$-adic norm.} around $x$. Construct a graph with vertices $B(x,r)$ and insert edges connecting $B(x,r)$ and $B(x',r')$ with $B(x,r) \subset B(x',r')$ and $[r,r'] \cap \abs{k^\times} = \set{r,r'}$. This graph is a simplicial tree, called the \defi{Bruhat-Tits-Tree} $\B(k)$. The ends\footnote{An \defi{end of a graph} is an infinite ray up to finitely many edges.} of this graph correspond bijectively to the points in $\P^1(\hat{k})$, where $\hat{k}$ denotes the completion of $k$. The action of $\PGL_2(k)$ on $\P^1(k)$ can be continued to an action on $\B(k)$, and we can modify $\B(k)$ by adding vertices such that this action is without inversion.

Let $\gamma \in \PGL_2(k)$ be hyperbolic or elliptic with two fixed points in $\P^1(k)$. In this case we define the \defi{axis} $\A(\gamma)$ to be the infinite path connecting the two ends of $\B(k)$ corresponding to the fixed points of $\gamma$. 
A hyperbolic element $\gamma \in \PGL_2$ acts on $\A(\gamma)$ by shifting the whole axis towards the end corresponding to the attracting fixpoint of $\gamma$. An elliptic element fixes $\A(\gamma)$ pointwise. It has additional fixed points in $\B(k)$ if and only if $\ord(\gamma)$ is a power of $p$ (this is made more precise in \zitat{Lemma 3}{Her}).
\pagebreak

Let $G$ be a finitely generated discontinuous subgroup of $\PGL_2(\C_p)$ and let $F(\gamma)$ be the set of the two fixed points of $\gamma \in G$ in $\P^1(\C_p)$. Now let $k$ be the extension of $\Q_p$ generated by the coefficients of the generators of $G$ and the fixed points of representatives of every conjugacy class of elliptic elements in $G$. There are only finitely many such conjugacy classes, hence $k$ is a finitely generated field extension of $\Q_p$ and we can therefore construct the Bruhat-Tits-Tree $\B(k)$ as described above. Note that by construction $G \subset \PGL_2(k)$ and $F(\gamma) \subset \P^1(k)$ for all elliptic elements $\gamma \in G$. Moreover we have $F(\gamma) \subset \P^1(\hat{k})$ for every hyperbolic element $\gamma \in G$ because the endpoints of $\A(\gamma)$ correspond to the fixed points of $\gamma$. 

As $G$ is discontinuous, $G$ contains only hyperbolic and elliptic elements (\zitat{Lemma 4.2}{Kato2}). 
The set of all fixed points of $G$
\[F(G) := \bigcup_{\gamma\in G} F(\gamma)\]
is a $G$-invariant subset of $\P^1(\hat{k})$, therefore $G$ also acts on the subtree $\T^*(G)$ of $\B(k)$ generated by the ends corresponding to $F(G)$.
We can now construct the quotient graph $\G^*(G) := \T^*(G)/G$.  Each axis of a hyperbolic element will be mapped to a circle in $\G^*(G)$, while each end of $\T^*(G)$ corresponding to a fixed point of an elliptic element but not to a fixed point of a hyperbolic element will be mapped to an end of $\G^*(G)$. We now turn $\G^*(G)$ into a graph of groups\footnote{In this chapter a \emph{graph} will always mean a \emph{graph of groups} to shorten notation. For the definition of a \defi{graph of groups} we refer to \zitat{I.4.4}{Ser}.} by labeling the image of a vertex resp. edge $x \in \T^*(G)$ with the conjugacy class of the stabilizing group $G_x$ of $x$. 

The graph $\G^*(G)$ contains a lot of useful information about the Mumford curve $\Omega(G)/G$. The ramification points of the covering $\Omega(G) \to \Omega(G)/G$ are the fixed points of elliptic elements of $G$  which are not fixed points of hyperbolic elements (\zitat{Prop. 5.6.2}{Kato2}). Therefore the branch points correspond bijectively to the ends of $\G^*(G)$. The stabilizing group of such an end is a cyclic group whose order equals the corresponding ramification index. And by studying the action of the hyperbolic elements, we find that the genus\footnote{By the \defi{genus of a graph} we mean its first Betti number.} of $\G^*(G)$ equals the genus of $\Omega(G)/G$ (\zitat{\S5.6.0}{Kato2}).

\begin{defn}
We define a \defindex{$p$-adic origami}{p-adic origami@$p$-adic origami} to be a covering of Mumford curves $X \to E$ ramified above at most one point with $g(E)=1$.
\end{defn}

Starting from a Mumford curve  $X= \Omega/\Gamma$ for a Schottky group $\Gamma$ we will later consider the quotient map to $E=\Omega/G$ for an extension $G$ of $\Gamma$. As $\Gamma$ is free the map $\Omega \to \Omega/\Gamma$ is unramified, therefore the branch points of $X \to E$ are equal to those of $\Omega \to \Omega/G$. Thus both necessary informations for $X \to E$ to be a $p$-adic origami (the genus of $E$ and the number of branch points) are coded in the quotient graph $\G^*(G)$ (as its genus and the number of its ends). We will now give two examples of how we can use this to construct $p$-adic origamis.

\pagebreak
\begin{ex}\label{ex:dihedral_origami}
Let $p > 5$ and $n \in \N$ be odd, $\zeta \in \C_p$ be a primitive $n$-th root of unity, $q \in \C_p$ with $\abs{q} < \abs{1-\zeta}$ and set
\[
\delta = \mattwo{\zeta}001 \quad 
\sigma = \mattwo0110 \quad 
\gamma = \mattwo{1+q}{1-q}{1-q}{1+q}.
\]
Thus $\delta$ is elliptic of order $n$ with fixed points $0$ and $\infty$, the involution $\sigma$ exchanges the fixed points of $\delta$ and has fixed points $1$ and $-1$, and $\gamma$ is hyperbolic with the same fixed points as $\sigma$. Then we have
 \vspace{-3pt}
\begin{enumerate}[$\bullet$] 
\setlength{\itemsep}{1pt}
  \setlength{\parskip}{0pt}
  \setlength{\parsep}{0pt}
\item $\gamma\sigma=\sigma\gamma$ and $\delta\sigma=\sigma\delta^{-1}$.
\item $\Delta := \gen{\delta,\sigma}$ is the dihedral group $D_n$ and fixes a single vertex $\A(\delta)\cap\A(\sigma)$.
\item $\Gamma := \gen{\delta^i\gamma\delta^{-i}: i\in\set{0,\dots,n-1}}$ is a Schottky group on $n$ free generators.
\item $\Gamma$ is a normal subgroup of $G := \gen{\delta,\sigma,\gamma}$ of index $2n$, hence $\Omega(G)=\Omega(\Gamma)=:\Omega$. It is the kernel of the map $\varphi : G \to \Delta$ defined by $\varphi|_\Delta = \id$ and $\varphi(\gamma)=1$. 
\item The quotient graph $\G^*$ of $\Omega/G$ is
\vspace{3mm}\begin{center}
\GraphOrigami{\Delta}{\left<\sigma\right>}{\left<\delta\right>}
\end{center}\vspace{-3mm}
where we use the arrow to indicate an end of this graph.
\item Since $\G^*$ has genus 1 and one end the map $\Omega/\Gamma \to \Omega/G$ is a normal $p$-adic origami with Galois group $G/\Gamma \cong D_n$.
\end{enumerate}  \vspace{-3pt}
A more detailed investigation of this origami can be found in \zitat{Bemerkung 4.3}{Kre}.
\end{ex}

\begin{ex}\label{ex:tetrahedral_origami}
Let $p > 5$ and $\zeta \in \C_p$ be a third root of unity, $q \in \C_p$ with $\abs{q}$ small enough\footnote{This is made more precise in \zitat{Bemerkung 4.4}{Kre}.} and set
\[
\delta = \mattwo{\zeta}001 \;\;
\sigma = \mattwo{-1}121 \;\;
\gamma = \mattwo{q}001.
\]
Thus $\delta$ is elliptic of order $3$ with fixed points $0$ and $\infty$, and $\gamma$ is hyperbolic with the same fixed points. 
The fixed points of the involution $\sigma$ are $-\frac12(1\pm\sqrt{3})$, those of $\sigma\delta\sigma$ are $\sigma(0)=1$ and $\sigma(\infty)=-\frac12$. 
Then we have \vspace{-3pt}
\begin{enumerate}[$\bullet$]
\setlength{\itemsep}{1pt}
  \setlength{\parskip}{0pt}
  \setlength{\parsep}{0pt}
\item $\gamma\delta=\delta\gamma$ and $(\delta\sigma)^3=\id$.
\item $\Delta := \gen{\delta,\sigma}$ is the tetrahedral group $A_4$ and fixes a single vertex $\A(\delta)\cap\A(\sigma)$.
\item $\Gamma := \gen{\alpha\gamma\alpha^{-1}: \alpha \in T}$  is a Schottky group on $4$ free generators.
\item $\Gamma$ is a normal subgroup of $G := \gen{\delta,\sigma,\gamma}$ of index $12$, hence $\Omega(G)=\Omega(\Gamma)=:\Omega$. It is the kernel of the map $\varphi : G \to \Delta$ defined by $\varphi|_\Delta = \id$ and $\varphi(\gamma)=1$. 
\item The quotient graph $\G^*$ of $\Omega/G$ is
\vspace{3mm}\begin{center}
\GraphOrigami{\Delta}{\left<\delta\right>}{\left<\sigma\right>}
\end{center}\vspace{-3mm}
\item Since $\G^*$ has genus 1 and one end the map $\Omega(\Gamma)/\Gamma \to \Omega(G)/G$ is a normal $p$-adic origami with Galois group $G/\Gamma \cong A_4$.
\end{enumerate}  \vspace{-3pt} \enlargethispage{12pt}
A more detailed investigation of this origami can be found in \zitat{Bemerkung 4.4}{Kre}.
\end{ex}

We will see in Section \ref{sct:p-adic_origamis} that the quotient graphs $\G^*$ of all non-trivial normal $p$-adic origamis look similar. We will then use this to investigate how the groups $G$ and $\Gamma$ have to be chosen such that the map $\Omega/\Gamma \to \Omega/G$ becomes a $p$-adic origami. But first we have to study the quotient graph $\G^*$ more closely.

\sect{Properties of the quotient graph}

A graph of groups $\G^*$ is called \defindex{$p$-realizable}{p-realizable@$p$-realizable}\footnote{In this section we will often just write $\defi{realizable}$ if the statements hold for arbitrary $p$.}, if there exists a finitely generated discontinuous group $G \subset \PGL_2(\C_p)$ with $\G^*=\G^*(G)$.
Let $\G\nt$ resp. $\G\nc$ be the subgraph of $\G^*$ containing only vertices and edges with non-trivial resp. non-cyclic groups.

\begin{thm}\label{thm:nb_endpoints}
The number of ends of a realizable graph $\G^*$ is 
\[n = \chi(\G\nc) + 2\chi(\G\nt)\]
where $\chi(\G)$ is the Euler-characteristic\footnote{Recall that the \defi{Euler-characteristic} (number of vertices minus number of edges) equals the difference of the first two Betti-numbers (number of connected components minus genus).} of a graph $\G$ (for infinite $\G$ we take the limit of $\chi$ for all finite subgraphs of $\G$).
\end{thm}
\begin{myproof}
Let $D$ resp. $d$ be the the number of vertices resp. edges in $\G\nc$. Then $\chi(\G\nc)=D-d$. Analogously let $C$ resp. $c$ be the number of vertices resp. edges in $\G\nt \setminus \G\nc$. Thus $\chi(\G\nt) = (C+D)-(c+d)$. Then we have to show
\[n = D-d+2((C+D)-(c+d)) = 2(C-c)+3(D-d)\]
Thus our statement is just a  reformulation of \zitat{Theorem 1}{Brad}.
\end{myproof}

\begin{lem}\label{lem:subgraphs}
Let $G$ be a finitely generated discontinuous group and let $\Ncal$ be a subgraph of $\G^*(G)$. Then there exists a subgroup $N$ of $G$ with quotient graph $\Ncal^* := \G^*(N) \supset \Ncal$, such that the difference between the two graphs $\Ncal$ and $\Ncal^*$ is contractible\footnote{An edge in a graph of groups may be \defindex{contracted}{contraction!of an edge in a graph} if it is not a loop and the inclusion of the edge group into one of its vertex groups is an isomorphism. After the contraction only the other vertex remains. Such a contraction does not change the fundamental group\footnotemark of the graph.}, except for the ends of $\Ncal^*$. 
\footnotetext{For the definition of the \defindex{fundamental group}{fundamental group!of a graph of groups} of a graph of groups we refer to \zitat{I.5.1}{Ser}.}
\end{lem}
\begin{myproof}
Choose a spanning tree of $\Ncal$ by deleting edges $\set{e_1,\dots,e_g}$ and let $\hat\Ncal$ be a preimage of this spanning tree in $\T^*(G)$. For each edge $e_i$ connecting vertices $v_i$ and $w_i$ let $\hat{e}_i$ be the lift of $e_i$ with $\hat{v}_i  \in \hat\Ncal$ and $\hat{e}'_i$ the lift with $\hat{w}'_i  \in \hat\Ncal$. The other endpoints $\hat{w}_i$ and $\hat{v}'_i$ of $\hat{e}_i$ resp. $\hat{e}'_i$ cannot be contained in $\hat{\Ncal}$ because otherwise $\hat{\Ncal}$ would contain a circle. Let $N$ be the subgroup of $G$ generated by all stabilizers of vertices in $\hat\Ncal$ and for each edge $e_i$ a hyperbolic element $\gamma_i$ mapping an $\hat{v}_i$ to $\hat{v}'_i$.

Thus $N$ is isomorphic to the fundamental group\footnotemark[\value{footnote}]  of the graph of groups $\Ncal$ by \zitat{I.5.4, Theorem 13}{Ser}. The stabilizers in $\Ncal$ do not change if we restrict the action from $G$ to $N$, neither do the identifications of vertices via the $\gamma_i$, hence the quotient graph $\T^*(G)/N$ contains $\Ncal$. Both graphs have the common fundamental group $N$, thus their difference cannot change their fundamental group, and therefore has to be contractible. 

The tree $\T^*(N)$ is contained in $\T^*(G)$, hence $\G^*(N)=\T^*(N)/N$ is contained in $T^*(G)/N$. Again both graphs have the common fundamental group $N$ and hence their difference is contractible.
\end{myproof}

\begin{prop}\label{prop:subgraphs}
Let $\G^*$ be a realizable graph, and let $\Ccal$ be a connected component of $\G\nc$. Then there exists a realizable graph $\Ncal^*$ with $\Ncal\nc = \Ccal$ (up to contractions) and $g(\Ncal\nc)=g(\Ncal^*)$.\vspace{-2pt}
\end{prop}
\begin{myproof}
Let $G$ be a finitely generated discontinuous group with $\G^*(G)=\G^*$. Subdivide all edges emanating from $\Ccal$ (which all have cyclic stabilizers), and let $\partial(\Ccal)$ be the set of all resulting edges in $\G^*\setminus \Ccal$ which still have a common vertex with $\Ccal$.
For the graph $\Ccal \cup \partial(\Ccal)$ Lemma \ref{lem:subgraphs} yields a graph $\Ncal^*$ with $\Ncal^* \supset \Ccal$. As $\Ncal^*$ and $\Ccal$ differ up to contraction only by ends, and ends are stabilized by cyclic groups, we get $\Ncal\nc=\Ccal$ (up to contractions) and $g(\Ncal\nc)=g(\Ccal)=g(\Ncal^*)$.\vspace{-2pt}
\end{myproof}

\begin{prop}\label{prop:circle_contains_cyclic_edge}
Let $\G$ be a connected graph of noncyclic groups which sa\-tis\-fies $g(\G)>0$. Then there exists no realizable graph $\G^*$ with $\G\nc = \G$ (up to contractions) and $g(\G\nc)=g(\G^*)$.\vspace{-2pt}
\end{prop}
\begin{myproof}
Assume there is a finitely generated discontinuous group $G$ such that $\G^* := \G^*(G)$ has the stated properties. 
$\G\nt$ is connected, because if $\sigma$ and $\tau$ are elements of stabilizers of two different connected components of $\G\nt$, then $\sigma\tau$ is hyperbolic and its axis contains the path $p$ between the axes of $\sigma$ and $\tau$. The image of $\A(\sigma\tau)$ in $\G^*$ is a circle which contains $p$ and hence an edge with trivial stabilizer. This would imply $g(\G^*)>g(\G\nc)$ contrary to the assumption.

Thus we have $g(\G\nt)= g(\G\nc)$ and both $\G\nt$ and $\G\nc$ are connected, hence $\chi(\G\nt)=\chi(\G\nc)$. Then Theorem \ref{thm:nb_endpoints} states $\chi(\G\nt)\geq 0$. But we have $\chi(\G\nt)<1$ because $\G\nt$ is connected and has positive genus. We thus get $\chi(\G\nt)=\chi(\G\nc)=0$ and therefore $g(\G\nt)=g(\G\nc)=1$. Hence by Theorem \ref{thm:nb_endpoints} the graph $\G^*$ has no ends. 

$G$ contains a normal subgroup $\Gamma$  of finite index which is a Schottky group. As $\G^*$ has no ends, the covering $\Omega(G)/\Gamma \to \Omega(G)/G$ is unramified. As $g(\Omega(G)/G)=1$, we conclude $g(\Omega(G)/\Gamma)=1$ by Riemann-Hurwitz.
Therefore $\Gamma$ is generated by a single hyperbolic element $\gamma$. All the elements in $G$ have the same axis as $\gamma$ (because otherwise there would be ramification points). Therefore every finite subgroup of $G$ is cyclic, which contradicts the assumption.\vspace{-2pt}
\end{myproof}

\begin{prop}\label{prop:genus_Gnc}
Let $\G^*$ be a realizable graph. Then $g(\G\nc)=0$.\vspace{-2pt}
\end{prop}
\begin{myproof}
For every connected component of $\G\nc$ this follows from Propositions \ref{prop:subgraphs} and \ref{prop:circle_contains_cyclic_edge}.\vspace{-2pt}
\end{myproof}

\begin{defn} 
Let $G\subset\PGL_2(\C_p)$ be a discontinuous group, $g(\Omega(G)/G) = 0$ and $\Omega(G) \to \Omega(G)/G$ ramified over exactly three points with ramification indices $n_1,n_2,n_3$. Then we call $G$ a ($p$-adic) \defi{triangle group} of type  $\Delta(n_1,n_2,n_3)$.\vspace{-2pt}
\end{defn}

The graph $\G^*(\Delta(n_1,n_2,n_3))$ is a tree with exactly three ends, corresponding to the three branch points.
Conversely if $G$ is a discontinuous group, $\G\nt$ is a tree, and $\G\nc$ is connected, then $\chi(\G\nt)=\chi(\G\nc)=1$, and $G$ is a triangle group  by Theorem \ref{thm:nb_endpoints}.\vspace{-2pt}\enlargethispage{12pt}

\begin{thm}\label{thm:contracted_graph}
Let $\G^*$ be a realizable graph and $\Ccal$ be a connected component of $\G\nc$. Then the fundamental group of $\Ccal$ is a triangle group $\Delta$. This means that $\Ccal$ can be replaced by a single vertex with vertex group $\Delta$ without changing the fundamental group of $\G^*$.
\end{thm}\pagebreak[4]
\begin{myproof}
Let $\Ncal^*$ be the realizable graph associated to $\Ccal$ by Proposition \ref{prop:subgraphs}.
By Proposition \ref{prop:genus_Gnc} this graph has genus zero, so by Theorem \ref{thm:nb_endpoints} it has three ends. Therefore the discontinuous group $\Delta$ with quotient graph $\Ncal^*$ is a triangle group.
\end{myproof}

Now we know that a $p$-realizable graph $\G$ with $\G\nc \neq \emptyset$ can be made up of vertices with $p$-adic triangle groups connected by edges with cyclic groups, it becomes vitally important to find all triangle groups which can occur. 
Fortunately for $p > 5$ those triangle groups are well-known:

\begin{thm}\label{thm:triangle_groups}
For every $p$ there exist the classical spherical triangle groups (i.e. those with $\frac1{n_1}+\frac1{n_2}+\frac1{n_3}>1$): the dihedral group $D_n=\Delta(2,2,n)$, and the symmetry groups of the platonic solids $A_4=\Delta(2,3,3)$, $S_4 = \Delta(2,3,4)$ and $A_5 = \Delta(2,3,5)$. For $p > 5$ there are no other $p$-adic triangle groups.
 \end{thm}
\begin{myproof}
Let $\Delta$ be one of the given groups. It is a finite subgroup of $\PGL_2(\C_p)$ and hence discontinuous. Its quotient graph $\G^*(\Delta)$ consists up to contraction of only one vertex (otherwise $\Delta$ would be a non-trivial amalgam or HNN-extension of smaller groups, hence would not be finite). This vertex has to be fixed by the whole group, which is non-cyclic, hence $\chi(\G\nc)=\chi(\G\nt)=1$. Therefore $\Delta$ is a triangle group by Theorem \ref{thm:nb_endpoints}.

Now let $\Delta$ be a triangle group for $p > 5$. We have $g(\G^*(\Delta))=0$, hence $\chi(\G\nc)\geq 1$ and $\chi(\G\nt)\geq 1$. By Theorem \ref{thm:nb_endpoints} we have then $\chi(\G\nc)=\chi(\G\nt)=1$, hence both graphs are connected. One can show that for $p > 5$ all edges of a realizable tree of groups can be contracted, which leaves a single vertex. The stabilizer of a vertex always is a finite subgroup of $\PGL_2(\C_p)$ and hence either cyclic or isomorphic to one of the groups stated (for a proof we refer to \zitat{Satz 2.7}{Kre}).
\end{myproof}

For $p \leq 5 $ there are additional non-spherical triangle groups. Bradley, Kato and Voskuil are currently working on their classification \cite{BKV}. A preliminary version and an idea of the proofs can be found in \cite{Kato1}.

\begin{ex}\label{ex:triangle}
 For $p=5$ an elliptic element $\delta$ of order 5 fixes not only its axis, but also all vertices contained in a small tube around this axis. Thus if we start with a vertex $v\in \B$ on the axis of $\delta$ fixed by a dihedral group $D_{5}$ generated by an element $\sigma$ of order 2 and $\delta$, then $\delta$ fixes also other vertices on the axis of $\sigma$. These vertices have the stabilizer $\left<\sigma,\delta\right>\cong D_5$ and we can find two elements $\tau$ and $\tau'$ of order 3, each with an axis through one of those vertices but not through $v$, such that the stabilizers of these two vertices under the action of $G := \gen{\sigma,\delta,\tau,\tau'}$ are $\gen{\sigma,\delta,\tau}\cong A_5$ and $\gen{\sigma,\delta,\tau'}\cong A_5$ respectively. If we do all this carefully we can get a discontinuous group whose quotient graph looks like this:
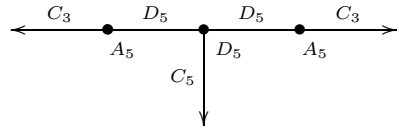
\begin{center}
\GraphTriangle
\end{center}
The generated group is thus a $p$-adic triangle group of type $\Delta(3,3,5)$. It is the fundamental group of the graph shown above, which is $A_5 *_{D_5} A_5$, where $*_{D_5}$ is the amalgamated product over the common subgroup $D_5$. Details about amalgams as fundamental groups of trees of groups can be found in \zitat{I.4.5}{Ser}.

One can generalize this example by starting with the dihedral group $D_{5n}$ for $n\in \N$. Then $\delta^n$ has order 5 and can be used to construct two stabilizers isomorphic to $A_5$. This results in a discontinuous group $A_5 *_{D_5} D_{5n} *_{D_5} A_5$, which is a $p$-adic triangle group of type $\Delta(3,3,5n)$.
\end{ex}

\sect{Normal $p$-adic origamis}\label{sct:p-adic_origamis}

After the preliminaries in the last two sections we are now ready to formulate our main result. We will restrict ourselves to ramified $p$-adic origamis, i.e. the case $g(X) > 1$. We are particularly interested in normal origamis, which we will now classify: 

\begin{thm}\label{thm:p-adic_origami}
Let $X \to E$ be a normal $p$-adic origami with $g(X)>1$. 
Then there is a discontinuous group $G$ and a Schottky group $\Gamma \NT G$ of finite index such that $X \cong \Omega/\Gamma$ and $E \cong \Omega/G$ with $\Omega := \Omega(\Gamma)=\Omega(G)$.

The group $G$ is isomorphic to the fundamental group of the graph of groups
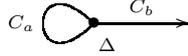
\vspace{1mm}
\begin{center}\vspace{4mm}
\GraphOrigami{\Delta}{C_a}{C_b}
\end{center}\vspace{-1mm}
where $\Delta$ is a $p$-adic triangle group of type $\Delta(a,a,b)$. This means that $G$ is isomorphic to the fundamental group of this graph.

Thus we get
\[G \cong \left<\Delta,\gamma;\, \gamma\alpha_1=\alpha_2\gamma\right> \text{ with $\alpha_i \in \Delta$ of order $a$.}\]
The Galois group of the origami is $G/\Gamma$. 
\end{thm}
\begin{myproof}
$X$ is a Mumford curve, thus there is a Schottky group $\Gamma \subset \PGL_2(\C_p)$ such that $X \cong \Omega(\Gamma)/\Gamma$. 
The automorphism group $\Aut X$ is isomorphic to $N/\Gamma$, where $N$ is the normalizer of $\Gamma$ in $\PGL_2(\C_p)$ (this is a theorem from \zitat{VII.2}{GvdP}). The Galois group of the covering $X\to E$ is a finite subgroup of $\Aut X$ and therefore takes the form $G/\Gamma$, where $\Gamma$ is a normal subgroup in $G \subseteq N$ of finite index. 
In this case $G$ is discontinuous and $\Omega(G)=\Omega(\Gamma)$.

The genus of $\G^*:=\G^*(G)$ equals the genus of $E$, which is 1.
The number of ends of $\G^*$ equals the number of branch points of the map $\Omega \to \Omega/G$. As the map $\Omega \to \Omega/\Gamma$ is unramified, this number equals the number of branch points of $X \cong \Omega/\Gamma \to E \cong \Omega/G$, which is also 1. Thus $\G^*$ is a realizable graph of genus one with one end. The stabilizer of this end is a cyclic group whose order equals the ramification index above the branch point.

We now prove that $\G\nc$ can be replaced by a vertex whose vertex group is a triangle group of the form $\Delta(a,a,b)$ ($a,b \in \N$): We know $g(\G\nt) \leq 1$, hence $\chi(\G\nt) \geq 0$, and the same holds for $\G\nc$. Theorem \ref{thm:nb_endpoints} states $\chi(\G\nc)+2\chi(\G\nt)=1$, therefore $\chi(\G\nc)=1$ and $\chi(\G\nt)=0$. As $g(\G\nt)\leq 1$ we conclude $g(\G\nt)=1$ and $\G\nt$ is connected. Prop. \ref{prop:genus_Gnc} states $g(\G\nc)=0$, therefore $\G\nc$ is connected as well. By Theorem \ref{thm:contracted_graph} we can replace $\G\nc$ by a single vertex $v$ whose vertex group is a triangle group $\Delta$. If we contract the rest of the graph as much as possible, we get an edge from $v$ to $v$ with a cyclic stabilizer. This stabilizer occurs therefore on two ends of $\G^*(\Delta)$ (this was $\Ncal^*$ in Theorem \ref{thm:contracted_graph}).
\end{myproof}

Note we have seen in Theorem \ref{thm:triangle_groups} that the spherical triangle groups of type $\Delta(a,a,b)$ are $D_n=\Delta(2,2,n)$ and $A_4=\Delta(2,3,3)$, and for $p>5$ there exist no other ones. For $p \leq 5$ there are additional possible triangle groups, for $p=5$ we have seen the type $\Delta(3,3,5n)$ in Example \ref{ex:triangle}.

We now know that normal $p$-adic origamis are always of the form $\Omega/\Gamma \to \Omega/G$, and we know quite well which groups $G$ can occur. It remains to investigate what groups $\Gamma$ are possible. The only restriction we have for $\Gamma$ is that it has to be a Schottky group of finite index and normal in $G$: As the covering $\Omega \to \Omega/\Gamma$ is always unramified the ramification of $\Omega/\Gamma \to \Omega/G$ is equal to the ramification of $\Omega \to \Omega/G$ and hence only depends on $G$. The genus of $\Omega/G$ also does not depend on the choice of $\Gamma$.

\begin{thm}\label{thm:schottky_equiv}
Let $G \subset \PGL_2(\C_p)$ be a finitely generated discontinuous group and $\Gamma$ be a normal subgroup of $G$ of finite index. Then the following statements are equivalent:
\begin{enumerate}[i)]
\setlength{\itemsep}{2pt}
  \setlength{\parskip}{0pt}
  \setlength{\parsep}{0pt}
\item $\Gamma$ is a Schottky group 
\item $\Gamma \cap G_i = \set{1}$ for every vertex group $G_i$ in $\G^*(G)$.
\end{enumerate}
\end{thm}
\begin{myproof}
i) $\RA$ ii) is easy: The vertex groups are finite, therefore every $g \in G_i$ has finite order. $\Gamma$  does not contain elements of finite order. Thus $\Gamma \cap G_i$ is trivial.
For ii) $\RA$ i) we proceed with three steps:

Step 1: Every element of $\Gamma$ has infinite order:
$G$ is the fundamental group of $\G^*(G)$, hence an HNN-extension of an amalgamated product of the $G_i$. If $g \in G$ has finite order $n > 1$, then $g$ is conjugated to a $g' \in G_i$ (see \zitat{IV.2.4 and IV.2.7}{LS}) with $\ord(g') = \ord(g) > 1$.  But by assumption $g' \not\in \Gamma$ and hence $g \not\in \Gamma$ as $\Gamma$ is normal in $G$.

Step 2: $\Gamma$ is free by Ihara's theorem (\zitat{I.1.5, Theorem 4}{Ser}):
$G$ acts on the tree $\T^*(G)$ with quotient graph $\G^*(G)$. For an $x \in \T^*(G)$ let $g \in \Stab_G(x)$ be non-trivial. Then $g$ has finite order, and with step 1 we see $g \not\in \Gamma$. Therefore the action on $\T^*(G)$ restricted to $\Gamma$ is free, thus $\Gamma$ is a free group by \zitat{\S 3.3}{Ser}.

Step 3: $\Gamma$ is a Schottky group:
$\Gamma$ is by definition a finite index normal subgroup of $G$. It is discontinuous because $G$ is and it contains no elements of finite order. It is finitely generated because $G$ is (by Reidemeister-Schreier, \zitat{Prop. II.4.2}{LS}). Thus we know that $\Gamma$ is a Schottky group. We can even find a finite set of free generators of $\Gamma$ by looking at its action on $\T^*(G)$: this action has a finite fundamental domain and this domain therefore has only finitely many neighboring translates. The set of these neighboring translates corresponds to a finite set of free generators for $\Gamma$.  
\end{myproof}

We are especially interested in the resulting Galois group $H := G/\Gamma$. Thus we now answer the question what choices of $\Gamma$ are possible if we fix this Galois group:
\pagebreak

\begin{cor}\label{cor:schottky_equiv}
Let $H$ be a finite group and $G \subset \PGL_2(\C_p)$ be a finitely generated discontinuous group. Further let $\Gamma$ be the kernel of a homomorphism $\varphi : G \to H$. Then the following statements are equivalent:
\begin{enumerate}[i)]
\setlength{\itemsep}{2pt}
  \setlength{\parskip}{0pt}
  \setlength{\parsep}{0pt}
\item $\Gamma$ is a Schottky group 
\item $\varphi\restrict{G_i}$ is injective for every vertex group $G_i$ in $\G^*(G)$.
\end{enumerate}
\end{cor}
\begin{myproof}
Follows with $\ker(\varphi\restrict{G_i}) = \Gamma \cap G_i$ from the Theorem.
\end{myproof}

\begin{ex}\label{ex:p-adic_origamis}

\begin{enumerate}[a)]
\setlength{\itemsep}{2pt}
  \setlength{\parskip}{0pt}
  \setlength{\parsep}{0pt}
\item Set $G_n := \left<D_n,\gamma;\, \gamma\sigma=\sigma\gamma\right>$ as in Example \ref{ex:dihedral_origami} (with $n$ odd and $\ord(\sigma)=2$). We extend this example: Choose $m \in \N$ and define $\varphi: G_n \to D_n \times C_m$ by $\varphi|_{D_n} = (\id,1)$ and $\varphi(\gamma) = (1,c)$ where $c$ is a generator of $C_m$. Then Corollary \ref{cor:schottky_equiv} states that $\Gamma' := \ker(\varphi)$ is a Schottky group and the Galois group of the $p$-adic origami $\Omega/\Gamma' \to \Omega/G$ is $D_n \times C_m$. Note that $\Gamma' \subseteq \Gamma := \ker(G_n \to D_n)$ for every $m$, thus we have a covering of origamis $\Omega/\Gamma' \to \Omega/\Gamma \to \Omega/G$. 
\item  Set $G := \left<A_4,\gamma;\, \gamma\delta=\delta\gamma\right>$ as in Example \ref{ex:tetrahedral_origami} (with $\ord(\delta)=3$). We extend this example as in a): Choose $m \in \N$ and define $\varphi: G_n \to A_4 \times C_m$ by $\varphi|_{A_4} = (\id,1)$ and $\varphi(\gamma) = (1,c)$ where $c$ is a generator of $C_m$. Then Corollary \ref{cor:schottky_equiv} states that $\Gamma' := \ker(\varphi)$ is a Schottky group and the Galois group of the $p$-adic origami $\Omega/\Gamma' \to \Omega/G$ is $A_4 \times C_m$. Note that $\Gamma' \subseteq \Gamma := \ker(G_n \to A_4)$ for every $m$, thus we have again a covering of origamis $\Omega/\Gamma' \to \Omega/\Gamma \to \Omega/G$. 
\item In Example \ref{ex:triangle} we have constructed the $5$-adic triangle group $\Delta(3,3,5)= A_5 *_{D_5} A_5$. The group
$G := \left<A_5 *_{D_5} A_5,\gamma;\, \gamma\delta_1=\delta_2\gamma\right>$, where the $\delta_i$ of order 3 are chosen out of the two different $A_5$-components, can be embedded into $\PGL_2(\bar\Q_5)$ (to show this one can use \zitat{Theorem II}{Kato2}). Then we can define $\varphi: G \to A_5$ as the identity on both $A_5$-components of the amalgamated product, and $\varphi(\gamma)=1$. This leads to a $5$-adic origami with Galois group $A_5$.
\item Take the group $G_5$ from part a) and consider the homomorphism $\varphi:G_5 \to \PSL_2(\F_{11})$ defined by
\[\varphi(\sigma)=\mattwo6715,\quad \varphi(\delta)=\mattwo0963, \quad \varphi(\gamma)=\mattwo3788\]
We can calculate that this is indeed a homomorphism by checking the relations of $G_5$ for the given images of $\varphi$. We can also check that $\varphi$ is surjective and that $\varphi|_{D_5}$ is injective. Hence $\ker(\varphi)$ is a Schottky group and the corresponding origami has Galois group $\PSL_2(\F_{11})$.
\end{enumerate}
\end{ex}

\sect{Automorphisms of $p$-adic origamis}\enlargethispage{12pt}

The Galois group of the covering $X \to E$ is a subgroup of the automorphism group $\Aut X$. For normal complex origamis we know from Section \ref{sct:complex} that the Galois group consists precisely of all possible translations. If the automorphism group is strictly larger than the Galois group, then there have to be automorphisms which are not translations, i.e. there is an automorphism which does not induce the identity but an involution on $E$. We now investigate the implications if this happens in the $p$-adic setting.

\begin{thm}\label{thm:p-adic_origami_aut}
In the situation of Theorem \ref{thm:p-adic_origami} let $\Aut(X)$ contain an element $\sigma$ of order 2 which induces a non-trivial automorphism $\bar{\sigma}$ of $E$ fixing the branch point of $X \to E$. Then there is a discontinuous group $H$ containing $G$ as normal subgroup of index 2, which is isomorphic to the fundamental group of the graph of groups 
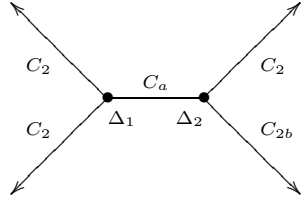
\begin{center}
\GraphOrigamiAut{\Delta_1}{\Delta_2}{C_a}{C_2}{C_2}{C_2}{C_{2b}}
\end{center}
where $\Delta_1$ is the $p$-adic triangle group of type $\Delta(2,2,a)$, i.e. $\Delta_1 \cong D_a$, and $\Delta_2$ is a $p$-adic triangle group of type $\Delta(2,a,2b)$ containing $\Delta$ of index 2.
\end{thm}
\begin{myproof}
Let $L$ be the subgroup of $\Aut X$ generated by $\sigma$ and the Galois group $\Gal(X/E)$. Every $\ell \in L\setminus \Gal(X/E)$ induces $\bar{\sigma}$ on $E$, thus $\ell \circ \sigma \in \Gal(X/E)$. Hence $L$ contains $\Gal(X/E)$ with index 2, and therefore as a normal subgroup.  As in the proof of Theorem \ref{thm:p-adic_origami} we have $L \cong H/\Gamma$ and for a discontinuous group $H$ and $\Gal(X/E) \cong G/\Gamma$ for a normal subgroup $G$ of index 2 in $H$ with $\Omega(H)=\Omega(G)=\Omega(\Gamma)$. Now $\Omega/H \cong E/\left<\sigma\right> =: P$.

The branch point of $X\to E$ is a fixed point of $\bar{\sigma}$ and therefore a ramification point of $E \to P$. By Riemann-Hurwitz this means that $g(P)=0$ and there are four ramification points of $E \to P$. As the degree of the map $E \to P$ is 2, and this is also the ramification index of the four ramification points, we know that there are also exactly four branch points. The composition $X\to P$ of the maps thus has four branch points. Over three of them the map $X \to E$ is unramified, therefore the corresponding ramification indices of $X \to P$ are $2$. Over the fourth branch point the map $X \to E$ is ramified with ramification index $b$, thus the total ramification index is $2b$. 

Now let $\Hcal^*$ be the quotient graph corresponding to $H$. Since $\Omega/H \cong P$ the graph $\Hcal^*$ is a realizable graph of genus zero with four ends. The stabilizer of one end is a cyclic group of order $2b$, the stabilizers of the other three ends are cyclic groups of order $2$.

We have $g(\Hcal^*)=0$, thus this also holds for all subgraphs of $\Hcal^*$. Therefore any Euler characteristic equals the number of connected components and we have $\chi(\Hcal\nc)\geq1$ and $\chi(\Hcal\nt)\geq 1$. As $\chi(\Hcal\nc)+2\chi(\Hcal\nt)=4$ by Theorem \ref{thm:nb_endpoints}, we conclude $\chi(\Hcal\nc)=2$ and $\chi(\Hcal\nt)=1$. Thus $\Hcal\nc$ has two connected components, which we can  by Theorem \ref{thm:contracted_graph} both replace by single vertices whose vertex groups are triangle groups $\Delta_1$ and $\Delta_2$. Furthermore $\Hcal\nt$ has to be connected, therefore those two triangle groups have to be connected by a path with a nontrivial cyclic stabilizer. This path can be contracted to a single edge.

It remains to find the connection between the stabilizing groups of $\G^*$ and $\Hcal^*$. We get the graph $\Hcal^*$ as the quotient of $\G^*$ by $\sigma$ (if we ignore the ends of both graphs). The single edge from $\G^*$ has to be mapped to itself and inverted (because otherwise there would still be a closed edge in $\Hcal^*$). Thus we have to insert a vertex on this edge and for constructing the quotient $\G^*$ we have to take one of the half edges. Therefore the stabilizer of this edge (as it is not fixed by $\sigma$) is the same as before (namely $C_a$) and the original vertex is fixed by $\Delta$ and $\sigma$.
\end{myproof}

\begin{ex}\label{ex:p-adic_aut}
Let $\zeta$ be a primitive 10-th root of unity and choose $\delta,\sigma$ and $\gamma$ as in Example \ref{ex:dihedral_origami}. Then the group $G := \gen{\delta^2,\sigma,\gamma}$ corresponds to the case $n=5$ of this example. But if we work out the quotient graph of $H := \gen{\delta,\sigma,\gamma}$ we note that while there still is a vertex with stabilizer $\gen{\sigma,\delta}\cong D_{10}$ now the element $\gamma\delta^5$ is elliptic of order 2, but does not fix this vertex. Instead there is now another vertex fixed by $\gamma\delta^5$ on the axis of $\sigma$, its stabilizer is therefore $\gen{\sigma,\gamma\delta^5}\cong D_2$. This means that $\gamma$ does not create a circle in the graph any more, but the quotient graph becomes
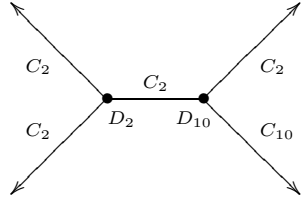
\begin{center}
\GraphOrigamiAut{D_2}{D_{10}}{C_2}{C_2}{C_2}{C_2}{C_{10}}
\end{center}
Now we define a homomorphism $\varphi:H \to \PSL_2(\F_{11}) \times \Z/2\Z$ by 
\[\varphi(\sigma)=\left(\mattwo6715,0\right),\quad \varphi(\delta)=\left(\mattwo864{10},1\right), \quad \varphi(\gamma)=\left(\mattwo3788,0\right).\]
As in Example \ref{ex:p-adic_origamis} d) we can check that this is indeed a homomorphism and is injective when restricted to the vertex groups. Note that $\varphi(\delta^2)=\left(\smat{0 & 9 \\ 6 & 3},0\right)$, thus $\varphi|_G$ is exactly the homomorphism considered in Example \ref{ex:p-adic_origamis} d), and we have $\ker(\varphi) = \ker(\varphi|_G) = \Gamma$. Thus the  $p$-adic origami $\Omega/\Gamma \to \Omega/G$ can be extended to $\Omega/\Gamma \to \Omega/H$ with Galois group $H/\Gamma \cong \PSL_2(\F_{11})\times \Z/2\Z$. This means that the automorphism group of this origami contains the group $\PSL_2(\F_{11})\times \Z/2\Z$.
\end{ex}


\sect{Complex and $p$-adic origamis}
Now we want to connect the complex and the p-adic worlds: An origami-curve in $\Mcal_{g,\C}$ is always defined over $\Qbar$, and thus defines also a curve in $\Mcal_{g,\Qbar}$, and this curve can in turn be interpreted as  a curve in $\Mcal_{g,\C_p}$. Now we ask the question: Does this curve intersect the subspace of $\Mcal_{g,\C_p}$ containing Mumford curves? 

The resulting points in $\Mcal_{g,\C_p}$ are still curves which cover an elliptic curve with only one branch point. Thus those curves are Mumford curves if and only if they occur as $p$-adic origamis. We have introduced several invariants of origami-curves in \cite{KreDiss}, some of which turn up in both worlds: The ramification indices, the Galois group of a normal origami and its automorphism group. By the Lefschetz principle (\zitat{Appendix}{Lef}) these algebraic properties of the origamis coincide over $\C$ and over $\C_p$. In some cases this is enough information to identify the complex origami-curve which belongs to a given $p$-adic origami.

Let $X_\C \to E_\C$ be an origami over $\C$. 
We can write $E_\C = E_{\C,\tau} = \C/(\Z+\tau\Z)$ with $\tau \in \H$, where $0$ is the only branch point of $X_\C \to E_\C$.
We have the Weierstrass-covering $\wp: E_{\C,\tau} \to \P^1(\C)$. For any ramification point $x\neq 0$ we have $\wp'(x)^2 = 4\wp^3(x)-g_2\wp(x)-g_3 = 0$, so if $g_2, g_3 \in \Qbar$ we get $\wp(x) \in \Qbar$. As $\wp(0)=\infty$ Belyi's theorem would then imply that both $X_\C$ and $E_\C$ are defined over $\Qbar$. Therefore we then have $X_\Qbar$ and $E_\Qbar$ over $\Qbar$ with the following diagram of base changes:
\[\begin{CD}
X_\Qbar \times_\Qbar \C = X_\C @>>> X_\Qbar @<<< X_{\C_p} = X_\Qbar \times_\Qbar \C_p\\
@VVV @VVV @VVV \\
E_\Qbar \times_\Qbar \C = E_\C @>>> E_\Qbar @<<< E_{\C_p} = E_\Qbar \times_\Qbar \C_p\\
\end{CD}\]

By varying $\tau \in \H$ we get a curve in the moduli-space $\Mcal_{g,\C}$, which leads to a curve in $\Mcal_{g,\Qbar}$, which itself can be considered as a subset of a curve in $\Mcal_{g,\C_p}$. This curve may or may not intersect the subset of Mumford curves in $\Mcal_{g,\C_p}$.

\begin{ex}\label{ex:con:kappes1}
Consider the following origami:
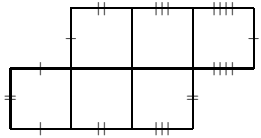
\begin{center}
\KappesOrigami
\end{center}

Kappes has proven in \zitat{Theorem IV.3.7}{Kap} that the origami-curve in $\Mcal_{g,\C}$ of this origami contains all the curves birationally equivalent to
\[y^2=(x^2-1)(x^2-\lambda^2)\left(x^2-\left(\tfrac{\lambda}{\lambda+1}\right)^2\right)\]
for $\lambda \in \C \setminus \set{0,\pm 1,-\frac12,-2}$. If we now restrict the choice of $\lambda$ to $\Qbar$, we get a curve in $\Mcal_{g,\Qbar}$. We can now change the base of this curve to $\C_p$ for an arbitrary prime $p$. This will result in a curve in $\Mcal_{g,\C_p}$. Does this curve intersect the subspace of $\Mcal_{g,\C_p}$ containing Mumford curves?

Fortunately \zitat{Theorem 4.3}{Brad-Cyc} offers a criterion for a hyperelliptic curve $X$ to be a Mumford curve: This is the case if and only if the branch points of $X \to \P^1$ (in our case $\pm1, \pm\lambda$ and $\pm \frac{\lambda}{\lambda+1}$) can be matched into pairs $(a_i,b_i)$ such that $\P^1$ can be covered by annuli $U_i$ each containing exactly one of those pairs. In our case we consider only $p>2$, set $\lambda := q-1$ for any $q \in \C_p^\times$ with $\abs{q}<1$ and match the points as follows:
\begin{align*}
a_1 &= 1,& b_1 &= -\lambda =1-q & \RA & \abs{a_1-b_1} = \abs{q} < 1 \\
a_2 &= -1,& b_2 &= \lambda =q-1 & \RA & \abs{a_2-b_2} = \abs{q} < 1 \\
a_3 &= \tfrac{\lambda}{\lambda+1} = \tfrac1q(q-1), 
&b_3 &= -\tfrac1q(q-1) & \RA & \abs{a_3} = \abs{b_3}= \abs{\tfrac1q}> 1
\end{align*}
Thus we can choose $U_1 = B(1,1)\setminus B(-1,\abs{q})$, $U_2 = B(1,1)\setminus B(1,\abs{q})$ and $U_3 = \P^1\setminus B(1,1)$ to get the desired covering\footnote{The balls $B(x,r)$ were defined in Section \ref{sct:discontinuous}.}. 
\end{ex}
\pagebreak

In general it is almost impossible to find the equation of a given origami. Therefore we would like a simpler approach, and we will do this the other way round: We have already constructed $p$-adic origamis, now we try to match them to the corresponding complex origami-curve. We will do this by matching some of the corresponding invariants.

\sect{Galois groups with a unique curve}\label{sct:unique_curve}
Some Galois groups occur only for a single origami-curve over $\C$. We will now prove that this is the case for the Galois groups $D_n \times \Z/m\Z$ and $A_4 \times \Z/m\Z$, which occurred as Galois groups of $p$-adic origamis in Example \ref{ex:p-adic_origamis} a) and b). Recall from Section \ref{sct:complex} that for a normal origami with monodromy $f : F_2 \onto G$ the other origamis on the same curve arise from $f$ by concatenation of elements of $\Aut(F_2)$ and $\Aut(G)$.

\begin{lem}\label{lem:cyclic_origami_01}
Let $f : F_2 \to  \Z/m\Z$ be surjective. Then up to an automorphism of $F_2$ we can assume $f(x)=1$ and $f(y)=0$.
\end{lem}
\begin{myproof}
Let $c_x, c_y \in \N$ with $c_x= f(x)$ and $c_y = f(y)$ in $\Z/m\Z$.

We prove first that we can choose the representatives $c_x$ and $c_y$ coprime: Let $p_i$ be the prime factors of $c_x$ and set 
\[c_y' := c_y + m \cdot \prod_{p_i \nmid c_y}p_i \]
Assume that there is a $p_i$ which is a factor of $c_y'$. If $p_i \mid c_y$ then $p_i$ would also have to be a factor of the right-hand summand, and as it is not contained in the product we would then have $p_i \mid m$. But this contradicts $\gen{c_x,c_y} = \Z/m\Z$. If on the other hand $p_i \nmid c_y$, then $p_i$ would be a factor of the right-hand summand but not of the left one, which would contradict the assumption.

Therefore no $p_i$ is a factor of $c_y'$, thus $\gcd(c_x,c_y') = 1$. We can replace $c_y$ by $c_y'$ as both are equivalent modulo $m$.

Now $\gcd(c_x,c_y) = 1$, thus there exist $a,b \in \Z$ with $ac_x+bc_y = 1$. Set
\[A := \mattwo{a}{-c_y}{b}{c_x} \in \SL_2(\Z)\]
and let $\varphi \in \Aut(F_2)$ be a lift of $A$.
As $\Z/m\Z$ is abelian we get a commutative diagram
\begin{center} 
\commdiag{@!=2mm}{
      F_2 \ar[rr]^\varphi \ar[rd] && F_2 \ar[rr]^f \ar[rd]  &     &  \Z/m\Z  \\
             			 &\Z^2 \ar[rr]^{\bar{\varphi}} &  &  \Z^2  \ar[ur]_{\bar{f}} &
}
\end{center}
where $\bar{\varphi}$ is the multiplication with $A$.

Thus 
\begin{align*}
(f\circ\varphi)(x) &= (\bar{f}\circ\bar{\varphi})\left(\mat{1\\0}\right) = \bar{f}\left(\mat{a \\ b}\right) = ac_x+bc_y = 1 \\
(f\circ\varphi)(y) &= (\bar{f}\circ\bar{\varphi})\left(\mat{0\\1}\right) = \bar{f}\left(\mat{-c_y \\ c_x}\right) = -c_yc_x+c_xc_y = 0 
\end{align*}
\end{myproof}

\begin{prop}
Let $n,m \in \N$ and $n$ be odd. Up to an automorphism of $F_2$ there exists only one curve of origamis with Galois group $D_n \times \Z/m\Z$.
\end{prop}
\begin{myproof}
Choose $\sigma, \tau \in D_n$ with $\gen{\sigma,\tau}=D_n$ and $\ord(\sigma)=\ord(\tau)=2$, and set $\delta = \sigma\tau$.

Let $f : F_2 \to D_n \times \Z/m\Z$ be the monodromy of a normal origami. By Lemma \ref{lem:cyclic_origami_01} we can apply an automorphism of $F_2$ to get $f(x)=(\alpha,1)$ and $f(y)=(\beta,0)$ for some $\alpha,\beta \in D_n$. As $f$ has to be surjective we know $\gen{\alpha,\beta} = D_n$. Up to an automorphism of $D_n$ we have three cases:
\begin{enumerate}[i)]
\setlength{\itemsep}{2pt}
  \setlength{\parskip}{0pt}
  \setlength{\parsep}{0pt}
\item $\ord(\alpha)=2, \ord(\beta)=2$, w.l.o.g. $(\alpha,\beta)=(\tau,\sigma)$
\item $\ord(\alpha)=n, \ord(\beta)=2$, w.l.o.g. $(\alpha,\beta)=(\delta,\sigma)$
\item $\ord(\alpha)=2, \ord(\beta)=n$, w.l.o.g. $(\alpha,\beta)=(\tau,\delta)$
\end{enumerate}

We can apply $\varphi \in \Aut(F_2)$ with $x \mapsto yx$ and $y \mapsto y$ to get from i) to ii):
\begin{align*}
(f \circ \varphi)(x) &= f(yx) = (\sigma\tau,0+1) = (\delta,1) \text{ and } \\
(f \circ \varphi)(y) &= f(y) = (\sigma,0) 
\end{align*}

For odd $m$ we can apply $\varphi \in \Aut(F_2)$ with $x \mapsto x$ and $y \mapsto yx^m$ to get from i) to iii):
\begin{align*}
(f \circ \varphi)(x) &= f(x) = (\tau,1) \text{ and }  \\
(f \circ \varphi)(y) &= f(yx^m) = (\sigma\tau^m, 0+m) = (\delta,0)
\end{align*}

For even $m$ case iii) is not possible, as $f$ would not be surjective: Assume there is a preimage $z \in f^{-1}(\id,1)$. We know that $f(xy)=f(y^{-1}x)$ (because this holds in both components). Thus we can choose $z$ of the form $x^ay^b$, its image is $(\tau^a\delta^b,a)$. Now for $f(z)_1 = \id$ we need $a$ even, but for $f(z)_2 = 1$  we need $a=1$ in $\Z/m\Z$. For even $m$ this is impossible. Therefore  there is no preimage of $(\id,1)$ in this case.
\end{myproof}

\begin{prop}\label{prop:con:single_tetrahedral}
Up to an automorphism of $F_2$ there exists, for any $m \in \N$, only one curve of origamis with Galois group $A_4 \times \Z/m\Z$.
\end{prop}
\begin{myproof}
Let $f : F_2 \to A_4 \times \Z/m\Z$ be the monodromy of a normal origami. By Lemma \ref{lem:cyclic_origami_01} we can apply an automorphism of $F_2$ to get $f(x)=(\alpha,1)$ and $f(y)=(\beta,0)$. As $f$ has to be surjective we know $\gen{\alpha,\beta} = A_4$. Up to an isomorphism of $A_4$ we have four cases:
\pagebreak[2]
\begin{enumerate}[i)]
\setlength{\itemsep}{2pt}
  \setlength{\parskip}{0pt}
  \setlength{\parsep}{0pt}
\item $\ord(\alpha)=2, \ord(\beta)=3$, w.l.o.g. $(\alpha,\beta)=((1\,2)(3\,4),(2\,3\,4))$
\item $\ord(\alpha)=3, \ord(\beta)=2$, w.l.o.g. $(\alpha,\beta)=((1\,2\,3),(1\,2)(3\,4))$
\item $\ord(\alpha)= \ord(\beta)=3$, $\ord(\alpha\beta)=2$, w.l.o.g. $(\alpha,\beta)=((1\,2\,3),(2\,3\,4))$
\item $\ord(\alpha)= \ord(\beta)=3$, $\ord(\alpha\beta)=3$, w.l.o.g. $(\alpha,\beta)=((1\,2\,4),(2\,3\,4))$
\end{enumerate}

We can apply $\varphi \in \Aut(F_2)$ with $x \mapsto xy$ and $y \mapsto y$ to get from iii) to i):
\begin{align*}
(f \circ \varphi)(x) &= f(xy) = ((1\,2\,3)(2\,3\,4),1) = ((1\,2)(3\,4),1) \text{ and } \\
(f \circ \varphi)(y) &= f(y) = ((2\,3\,4),0) 
\end{align*}

We can apply $\varphi \in \Aut(F_2)$ with $x \mapsto xy^{-1}$ and $y \mapsto y$ to get from iii) to iv):
\begin{align*}
(f \circ \varphi)(x) &= f(xy^{-1}) = ((1\,2\,3)(2\,4\,3),1) = ((1\,2\,4),1) \text{ and } \\
(f \circ \varphi)(y) &= f(y) = ((2\,3\,4),0) 
\end{align*}

If $m$ is not divisible by 3 choose $k \in \N$ such that $km \equiv 1 \mod 3$. We can then apply $\varphi \in \Aut(F_2)$ with $x \mapsto x$ and $y \mapsto x^{km}y$ to get from iii) to ii):
\begin{align*}
(f \circ \varphi)(x) &= f(x) = ((1\,2\,3),1) \text{ and }  \\
(f \circ \varphi)(y) &= f(x^{km}y) = ((1\,2\,3)(2\,3\,4), 0+km) = ((1\,2)(3\,4),0)
\end{align*}

If $m$ is divisible by 3, case ii) is not possible, as $f$ would not be surjective: 
Assume there is a preimage $z \in f^{-1}(\id,1)$. As the first component of $f(z)$ is $\id \in A_4$ the presentation of $A_4$ tells us that the element $z$ is contained in the normal subgroup generated by $y^2,x^3$ and $(yx)^3$, and hence is a product of conjugates of those elements. But $f(y^2)=(\id,0)$ and $f(x^3)=f((yx)^3)=(\id,3)$. Therefore the second component of $f(z)$ is divisible by 3. As $m$ is also divisible by 3 this contradicts $f(z)=(\id,1)$.
\end{myproof}

\begin{ex}\label{ex:con:kappes2}
In Example \ref{ex:con:kappes1} we have shown that the origami $O$ defined by
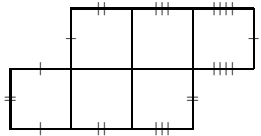
\begin{center}
\KappesOrigami
\end{center}
occurs as a $p$-adic origami. Now we have an alternative way of showing this: Let $O'$ be the minimal normal cover of $O$. The Galois group of $O'$ is $A_4$. Proposition \ref{prop:con:single_tetrahedral} tells us that the origami-curve of $O'$ is the only curve with this Galois group. 

In Example \ref{ex:tetrahedral_origami} we constructed a $p$-adic origami with Galois group $A_4$; in fact one for every suitably chosen $q \in \P^1(\C_p)$. If one of those origamis $X'$ is (as an algebraic curve) defined over $\Qbar$ and hence over $\C$ then it has to occur on the origami-curve of $O'$. The covering $O'\to O$ leads to a morphism $X' \to X$, where $X$ is on the origami-curve of $O$. As $X'$ is a Mumford curve the same holds for $X$ by \zitat{Satz 5.24}{Brad-Dip}.
\end{ex}

\sect{When the group is not enough}

In \zitat{2.2}{KreDiss} we have proposed an algorithm to find all complex origami-curves of normal origamis with a given Galois group $H$. This results in a set of representative origamis given as epimorphisms $f : F_2 \onto H$ as described in Section  \ref{sct:complex}. For most groups there is only one curve with the given Galois group, but there are cases where there are more than one.\footnote{There are 2386 groups with order less than or equal to 250 which can be generated by two elements. Of these there are only 30 where there is more than one curve with this Galois group, c.f. \zitat{A.3}{KreDiss}.} The smallest example for such a group is the group $A_5$ where there are two curves. Representatives are given by
\begin{align*}
f_1:F_2 \to A_5, \quad x&\mapsto (1\,5\,3\,4\,2), \quad y \mapsto (1\, 3\, 2\, 4\, 5)\quad \text {and}\\
f_2:F_2 \to A_5, \quad x&\mapsto (1\,5\,2\,4\,3), \quad y \mapsto (2\, 3)(4\, 5)
\end{align*}
We can easily see that those two origamis do not define the same curve: The ramification index of a normal origami is the order of $f(xyx^{-1}y^{-1})$. Hence the ramification indices of those two origamis are
\[
\ord(f_1(xyx^{-1}y^{-1})) = 5 \text { and } \ord(f_2(xyx^{-1}y^{-1})) = 3 
\]

\begin{ex}
In Example \ref{ex:p-adic_origamis} c) we have investigated a $5$-adic origami with Galois group $A_5$. 
Our $5$-adic origami had ramification index 5 (recall that this was the order of the cyclic stabilizer of the single end of the quotient graph), hence this corresponds to the curve of the origami defined by $f_1$.
\end{ex}

Sometimes fixing the ramification index makes the origami-curve unique. But there are still some cases where two curves have equal Galois groups and equal ramification index. An example for such a group is the group $\PSL_2(\F_7)$, where there are even four curves, represented by $f_i : F_2 \to \PSL_2(\F_7)$ with $f_i(x)=\sigma_i$ and $f_i(y)=\tau_i$ with
\begin{align*}
&\sigma_1=\sigma_2=\sigma_3 = \mattwo1254, \quad \sigma_4=\mattwo0616 \\
&\tau_1 = \mattwo6026, \quad \tau_2=\sigma_4, \quad \hspace{3pt}\tau_3= \mattwo1365, \quad \tau_4=\mattwo6334
\end{align*}
The two curves defined by $f_1$ and $f_2$ both have ramification index 4.

In \zitat{Prop. 3.13}{KreDiss} we show that the automorphism group of a normal origami with Galois group $H$ is either $H$ itself or isomorphic to $C_2 \ltimes_\Phi H$ where $\Phi: C_2 \to \Aut(H)$ maps the generator of the cyclic group $C_2$ to the automorphism $\varphi$ of $H$ defined by $\varphi(f(x))=f(x)^{-1}$ and $\varphi(f(y))=f(y)^{-1}$. In the example above with Galois group $\PSL_2(\F_7)$ and ramification index 4 the automorphism groups are isomorphic and hence can not be used to distinguish those two curves. But in some cases they are helpful:

\begin{ex}
In Example \ref{ex:p-adic_origamis} d) we have investigated a $5$-adic origami with Galois group $\PSL_2(\F_{11})$ with ramification index $5$. There are four possible origami-curves, represented by $f_i : F_2 \to \PSL_2(\F_{11})$ with $f_i(x)=\sigma_i$ and $f_i(y)=\tau_i$ with
\begin{align*}
&\sigma_1=\sigma_2=\sigma_3 = \mattwo4684,\quad \,\sigma_4=\mattwo{10}{10}76 \\
&\tau_1 = \sigma_4, \quad \tau_2 = \mattwo218{10},\quad  \tau_3= \mattwo1401, \quad \tau_4=\mattwo4{10}03
\end{align*}
The automorphism groups of $f_2$ and $f_4$ are isomorphic to $\PSL_2(\F_{11})\times \Z/2\Z$, while the other two are isomorphic to $\Aut(\PSL_2(\F_{11}))$.\
We have seen in Example \ref{ex:p-adic_aut} that the automorphism group of the $p$-adic origami contains $\PSL_2(\F_{11})\times \Z/2\Z$, thus the corresponding complex origami-curve is either defined by $f_2$ or $f_4$. We can presently not decide which of the two curves is the right one.
\end{ex}

\begin{acknowledgement}
 I want to express my gratitude to my thesis advisor Frank Herrlich for encouraging my work on Mumford curves and origamis and supporting me throughout the whole time. From Gabriela Schmith\"usen I have learned a lot about origamis. I also thank her for urging me to publish this paper. Patrick Erik Bradley introduced me to Kato's work on quotient graphs by handing me an early draft of \cite{BKV}, which led to the key ideas for the classification theorem. 
\end{acknowledgement}

\end{document}